\documentclass{article} 
\usepackage{graphicx}
\usepackage{natbib}

\usepackage{amsmath}        
\usepackage{amssymb}
\usepackage[english]{babel}
\usepackage{amsthm}
\usepackage{bbm}
\usepackage{url}
\usepackage{longtable}
\usepackage{booktabs}
\usepackage{caption}
\usepackage{algorithmic}
\usepackage{algorithm}
\usepackage{subfigure}
\usepackage{tabularx}

\newtheoremstyle{plainNoItalics}{}{}{\normalfont}{}{\bfseries}{.}{ }{}

\theoremstyle{plainNoItalics}
\newtheorem{theorem}{Theorem}

\newtheorem{problem}[theorem]{Problem}

\newtheorem*{theorem*}{Theorem}
\newtheorem*{proposition*}{Proposition}
\newtheorem*{lemma*}{Lemma}
\newtheorem*{corollary*}{Corollary}
\newtheorem*{remark*}{Remark}

\newtheorem*{observation*}{Observation}
\newtheorem*{example*}{Example}
\newtheorem*{assumption*}{Assumption}

\theoremstyle{definition}
\newtheorem*{definition*}{Definition}

\newenvironment{proofsketch}
  {\begin{proof}[Proof sketch]}
  {\end{proof}}

\newcommand{\R}[0]{R}

\begin{document}

\title
{Exploiting Packing Components \\ in General-Purpose Integer \\ Programming Solvers}

\author{Jakub Mare{\v c}ek \\ \footnotesize IBM Research -- Ireland, B3 F14 Damastown Campus, Dublin 15\\
\footnotesize The Republic of Ireland, \url{jakub.marecek@ie.ibm.com}
}

\maketitle

\begin{abstract}
The problem of packing boxes into a large box is often a part of a larger problem.
For example in furniture supply chain applications, one needs to decide what
trucks to use to transport furniture between production sites and 
distribution centers and stores, such that the furniture fits inside.
Such problems are often formulated and sometimes solved using general-purpose integer programming solvers.

This chapter studies the problem of identifying a compact formulation of the 
multi-dimensional packing component in a general instance of integer linear programming, reformulating it using the 
discretisation of Allen--Burke--Mare{\v c}ek, and and solving the extended reformulation.
Results on instances of up to 10000000 boxes are reported.
\end{abstract}

\section{Introduction}
\label{sec:Intro}

It is well known that one problem may have many integer linear programming formulations,
whose computational behaviour is widely different.
The problem of packing three-dimensional boxes into a larger box is a particularly striking
example. 
For the trivial problem of how many unit cubes can be packed into a $(k \times 1 \times 1)$ box, a widely known formulation of Chen/Padberg/Fasano makes it impossible to answer the question for as low a $k$ as $12$ within an hour using state of the art solvers (IBM ILOG CPLEX, FICO XPress MP, Gurobi Solver), despite the fact that (after pre-solve), the instance has only 616 rows and 253 columns and 4268 non-zeros. 
In contrast, a discretised (``space-indexed'') formulation of  \cite{Marecek2011} makes it possible to solve the instance with $k = 10000000$ within an hour, 
where the instance of linear programming had 10000002 rows, 10000000 columns, and 30000000 non-zeros.
This is due to the fact that the discretised linear programming relaxation provides a particularly strong bound.
On a large-scale benchmark of randomly generated instances, the difference between root linear programming relaxation value and optimum to optimum
has been 10.49 \% and 0.37 \% for the Chen/Padberg/Fasano and the formulation of \cite{Marecek2011}, respectively.

The Allen--Burke--Marecek formulation with adaptive discretisations is, however, often rather hard to formulate in an algebraic modelling language.
First, the choice of the discretisation of the larger box is hard in any case; in terms of computational complexity, finding the 
best possible discretisation is $\Delta_2^p$-Hard.
Second, algebraic modelling languages such as AMPL, GAMS, MOSEL, and OPL are ill-suited to the dynamic programming 
required by efficient discretisation algorithms.
Finally, non-uniform grids pose a major challenge in debugging the formulation and analysis of the solutions.
Within the modelling language, 
the discretisation is hence often chosen in an ad hoc manner, without considerations of optimality.

Instead, this chapter studies the related problems of 
identifying a packing component in the Chen/Padberg/Fasano formulation in 
a larger instance of integer linear programming,
automating the reformulation of Chen/Pad\-berg/Fasano formulation into the Allen--Burke--Marecek formulation, 
obtaining a reasonable discretisation in the process,
and translating the solutions thus obtained back to the original problem.
The main contributions are:

\begin{enumerate}
\item An overview of integer linear programming formulations of packing boxes into a larger box, covering both the formulation of Chen/Pad\-berg/Fasano and the discretisation of \cite{Marecek2011}.
\item A formal statement of the problems of extracting the component and row-block, respectively, which correspond to the formulation of Chen/Pad\-berg/Fasano for packing boxes into a larger box, from a general integer linear programming instance.
Surprisingly, we show that the extraction of the row-block is solvable by polynomial-time algorithms.
\item A formal statement of the problem of finding the best possible discretisation, i.e. reformulation of Chen/Padberg/Fasano to the formulation of \cite{Marecek2011}. 
We show that the problem is $\Delta_2^p$-Hard, but that there are very good
heuristics.
\item A novel computational study of the algorithms for the problems above.
There, we take the integer linear programming instance,
extract the row-block corresponding to the Chen/Padberg/Fasano formulation,
perform the discretisation, write out the  integer linear programming instance 
using the Allen--Burke--Marecek formulation, and solve it.
\end{enumerate}
A general-purpose integer linear programming solver using the above results
 could solve much larger instances of problems involving packing components
 than is currently possible, when they are formulated using the 
 Chen/Padberg/Fasano formulation.

\section{Background and Definitions}

In order to motivate the study of packing components, let us consider: 

\begin{problem}
\label{prob:container-loading}
The {\sc Precedence-Constrained Scheduling (PCS)}: Given integers $r, n \ge 1$, 
amounts $a_i > 1$ of resources $i = 1, 2, \ldots, r$ available,
resource requirements of $n$ jobs $D \in \R^{n \times r}$,
and $p$ pairs of numbers $P \subseteq \{ (i, j) \; | \; 1 \le i < j \le n \}$ expressing job $i$ should be executed prior to executing $j$,  
find the largest integer $k$ so that $k$ jobs can be executed using the resources available. 
\end{problem}

This problem on its own has numerous important applications, notably in extraction of natural resources (\cite{Bienstock2010,Moreno2010407}),
  where it is known as the open-pit mine production scheduling problem,
aerospace engineering (\cite{baldi2012three,fasano2013modeling}), 
 where one needs to balance the load of an aircraft, satellite, or similar,
and complex vehicle routing problems, 
 where both weight and volume of the load is considered (e.g., \cite{iori2007exact}).

Clearly, there is a packing component to Pre\-ce\-dence-Constrain\-ed Sche\-dul\-ing (Problem~\ref{prob:container-loading}).
Let us fix the order of six allowable rotations in dimension three arbitrarily and define:

\begin{problem}
The {\sc Container Loading Problem (CLP)}: Given dimensions of a large box (``container'') $x, y, z > 0$ and dimensions of
$n$ small boxes $D \in {\mathbbm{R}}^{n \times 3}$ with associated values $w \in {\mathbbm{R}}^{n}$, 
and specification of the allowed rotations $r = \{0, 1\}^{n \times 6}$, find the 
greatest $k \in \mathbbm{R}$ such that there is a packing of small boxes $I \subseteq \{1, 2, \ldots, n\}$ 
into the container 
with value $k = \sum_{i \in I} w_i$. 
The packed small boxes $I$ may be rotated in any of the allowed ways,
must not overlap, and no vertex can be outside of the container.
\end{problem}

\begin{problem}
The {\sc Van Loading Problem (VLP)}: Given dimensions of a large box (``van'') $x, y, z > 0$, maximum mass $p \ge 0$ it can hold (``payload''), 
dimensions of $n$ small boxes $D \in {\mathbbm{R}}^{n \times 3}$ with associated values $w \in {\mathbbm{R}}^{n}$, 
mass $m \in {\mathbbm{R}}^{n}$, 
and specification of the allowed rotations $r = \{0, 1\}^{n \times 6}$, find the greatest $k \in \mathbbm{R}$ 
such that there is a packing of small boxes $I \subseteq \{1, 2, \ldots, n\}$ into the container 
with value $k = \sum_{i \in I} w_i$ and mass $\sum_{i \in I} m_i \le p$.
The packed small boxes $I$ may be rotated in any of the allowed ways,
must not overlap, and no vertex can be outside of the container.
\end{problem}

Such problems are particularly challenging. 
Ever since the work of \cite{gilmore1965multistage}, there has been much research on extended formulations of 2D packing problems using the notion of patterns, e.g. \cite{madsen1979glass}. See \cite{ben2005cutting} for an excellent survey.
Only in the past decade or two has the attention focused to exact solvers for 3D packing problems (\cite{martello2000three}), where even the special case with rotations around combinations of axes in multiples of 90 degrees is NP-Hard to approximate (\cite{Chlebik09}).
Although there are a number of excellent heuristic solvers, 
the progress in exact solvers for the Container Loading Problem has been limited, so far.


\begin{table}[t]
	\caption{Notation used in this chapter, which matches \cite{Marecek2011}.}
	\centering
		\begin{tabular}{cl}
		Symbol & Meaning \\
		\midrule
			$n$ & The number of boxes. \\
			$H$ & A fixed axis, in the set $\{X,Y,Z\}$. \\
			$\alpha$ & An axis of a box, in the set $\{1,2,3\}$. \\
			$L_{\alpha i}$ & The length of axis $\alpha$ of box $i$. \\
			$l_{\alpha i}$ & The length of axis $\alpha$ of box $i$ halved. \\
			$D_H$ & The length of axis $H$ of the container. \\
			$w_i$ & The volume of box $i$ in the CLP. \\ 
     \bottomrule
		\end{tabular}
		\vskip 5mm
	\label{notationtable}
\end{table}

\paragraph{The Formulation of Chen/Padberg/Fasano}

The usual compact linear programming formulations provide only weak lower bounds.
\cite{Chen1995} introduced an integer linear programming formulation using the relative placement indicator:

\begin{align}
\lambda_{ij}^H & = \begin{cases}
   \; 1 & \text{if box } i \text{ precedes box } j \text{ along axis } H \\
   \; 0 & \text{otherwise} \notag 
 \end{cases}, \\
\delta_{\alpha i}^H & = \begin{cases}
   \; 1 & \text{if box } i \text{ is rotated so that axis } \alpha \text{ is parallel to fixed } H \\
   \; 0 & \text{otherwise} \notag 
 \end{cases}, \\ 
x_{i}^H & = \mbox{ absolute position of box } i \text{ along axis } H. \notag
\end{align}

Using the notation of Table~\ref{notationtable} and implicit quantification, it reads: 

\begin{flalign}
\max \; &\sum_{i=1}^n{\sum_H{w_i\delta_{1i}^H}}\\
\textrm{s.t. } \; 
& \sum_H{\delta_{2i}^H} = \sum_H{\delta_{1i}^H} \label{CP:eq6nzs:1} \\ 
& \sum_H{\delta_{1i}^H} = \sum_\alpha \delta_{\alpha i}^H \label{CP:eq6nzs:2} \\ 
&L_{1j(i)}\lambda^H_{j(i)i} + \sum_\alpha l_{\alpha i}\delta_{\alpha i}^H \leq x_i^H \label{CP:5nzs:1} \\ 
& x_i^H \leq \sum_\alpha{(D_H - l_{\alpha i})\delta_{\alpha i}^H} - L_{1j(i)}\lambda_{ij(i)}^H \label{CP:5nzs:2} \\ 
&D_H\lambda_{ji}^H + \sum_\alpha{l_{\alpha i}\delta_{\alpha i}^H} - \sum_\alpha{(D_H - l_{\alpha j})\delta_{\alpha j}^H} \leq x_i^H - x_j^H \label{CP:12nzs:1} \\ 
&x_i^H - x_j^H \leq \sum_\alpha{(D_H - l_{\alpha i})\delta_{\alpha i}^H} - \sum_\alpha{l_{\alpha j}\delta_{\alpha j}^H} - D_H\lambda_{ij}^H \label{CP:12nzs:2} \\ 
&\sum_H{(\lambda_{ij}^H + \lambda_{ji}^H)} \leq \sum_H{\delta_{1i}^H} \label{CP:9nzs:1} \\
&\sum_H{(\lambda_{ij}^H + \lambda_{ji}^H)} \leq \sum_H{\delta_{1j}^H} \label{CP:9nzs:2} \\ 
&\sum_H{\delta_{1i}^H} + \sum_H{\delta_{1j}^H} \leq 1 + \sum_H{(\lambda_{ij}^H + \lambda_{ji}^H)} \label{CP:12nzs} \\ 
&\sum_{i=1}^n{\sum_H{\left(\prod_\alpha {L_{\alpha i}}\right)\delta_{1i}^H \leq \prod_H{D_H}}} \label{CP:3nnzs} \\ 
&\delta_{\alpha i}^H \in \{0,1\}, \lambda_{ij}^H \in \{0,1\} \notag \\ 
& 
L_{1i}\leq L_{2i} \leq L_{3i}, j(i) \text{ such that } L_{1j(i)} = \max\{L_{1j}\} \text{ for } 1 \leq i \neq j \leq n. \notag
\end{flalign}

\begin{figure*}[p]
	\centering
	\caption{The $A$ matrix corresponding to the previously unsolved instance Pigeon-02 in the Chen/Padberg/Fasano formulation,
  with colour highlighting the absolute value of the coefficients.
  }
  \vskip 5mm
	\includegraphics[scale=1.1,trim={0cm 0cm 0cm 0cm},clip]{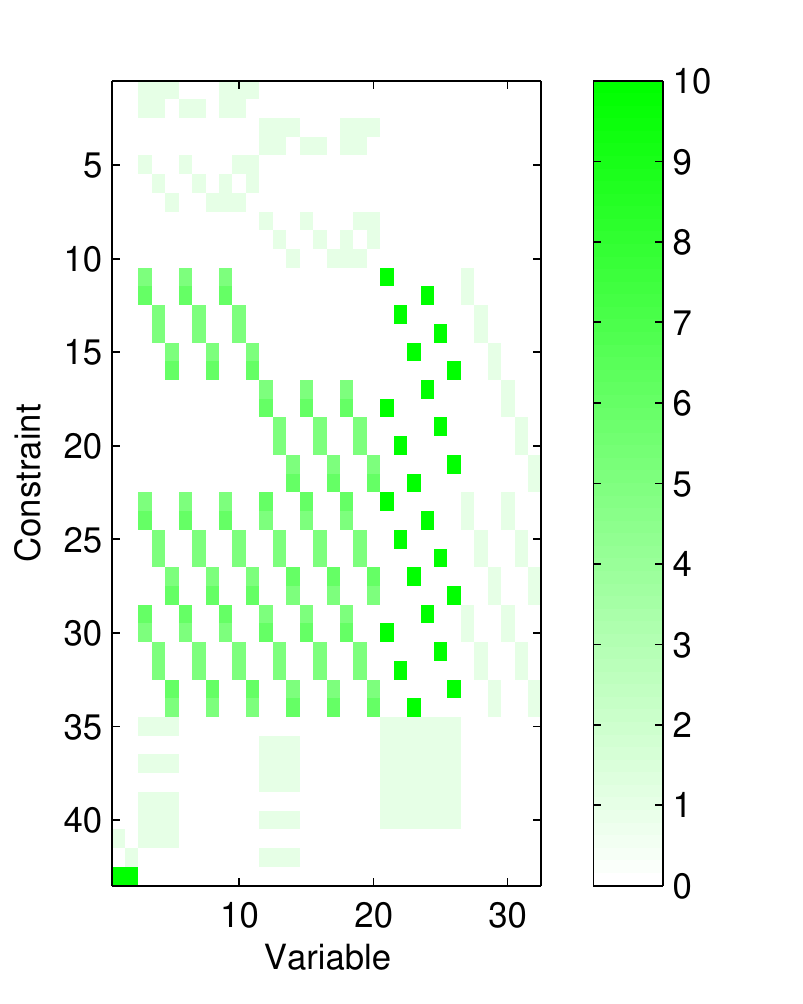}
	\label{pigeon02}
\end{figure*}

See Figure~\ref{pigeon02} for the sparsity pattern of a small instance, known as Pigeon-02,
where at most a single unit cube out of two can be packed into a single unit cube. 
The constraint matrix of Pigeon-02 is $43 \times 32$.
In the figure, the column ordering is given by placing $D$ first, $\delta$ second, $\lambda$ third, and $x$ at the end, with the highest-order-first indexing therein.
$D$ has dimension two, $\delta$ has dimension 18, $\lambda$ has dimension 6
and $x$ has dimension 6. 
In the figure, the row ordering is given by the order of constraints (\ref{CP:eq6nzs:1}--\ref{CP:3nnzs}) above.
Notice that a similar ordering of columns and rows is naturally produced by a parser of an algebraic modelling language, 
such as AMPL, GAMS, MOSEL, or OPL.


This formulation has been studied a number of times. 
Notably,
\cite{Fasano99,fasano2004mip,fasano2008mip}
suggested numerous improvements to the formulation. 
\cite{Padberg00} has studied properties of the formulation and, in particular, identified
the subsets of constraints with the integer property. 
\cite{Marecek2011} proposed further improvements, including symmetry-breaking constraints
and means of exploitation of properties of the rotations.
See \cite{Fasano2014a,Fasano2014b,Fasano2014c} for further extensions to Tetris-like items
and further references.

\begin{figure*}[p]
	\centering
	\caption{The $A$ matrix corresponding to instance Pigeon-12 in the Chen/Padberg/Fasano formulation,
  with colour highlighting the absolute value of the coefficients.
  }
  \vskip 5mm
	\includegraphics[scale=1.1,trim={0cm 0cm 0cm 0cm},clip]{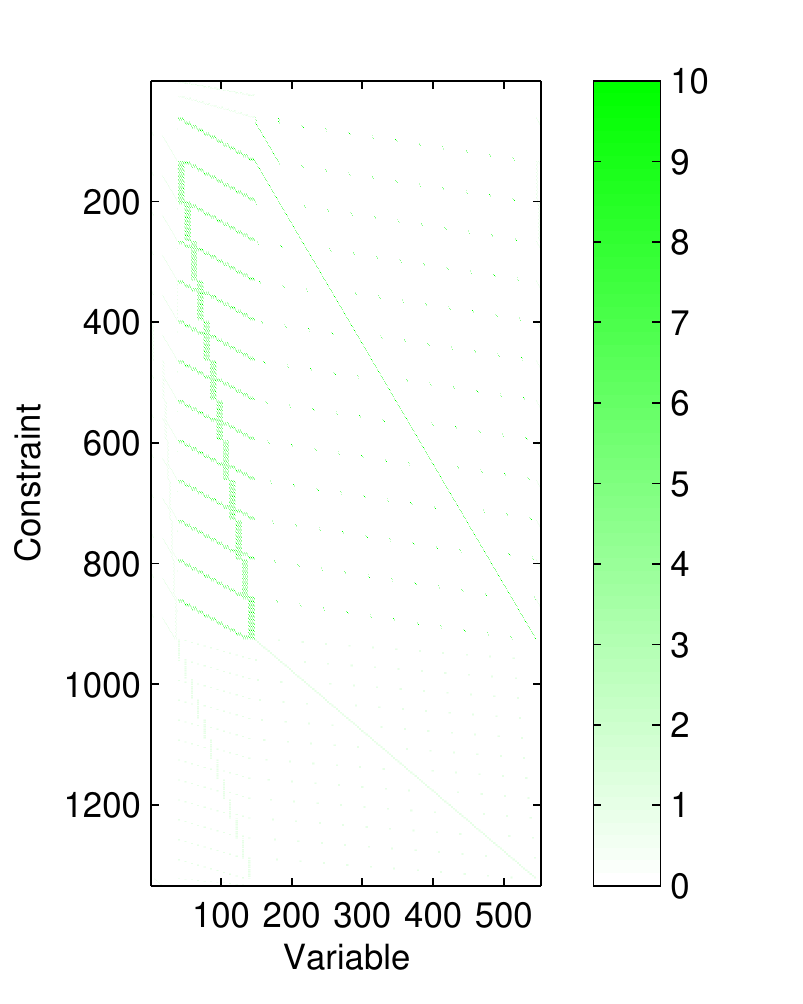}
	\label{pigeon12}
\end{figure*}

Nevertheless, whilst the addition of these constraints improves the performance somewhat, 
the formulation remains far from satisfactory.
As has been pointed out in Section~\ref{sec:Intro} and can be confirmed in Table~\ref{tab:pigeon-results}, 
Pigeon-$k$ becomes very challenging as the number $k$ of unit cubes to pack into $(k \times 1 \times 1)$  box grows.
See Figure~\ref{pigeon12} for the sparsity patten of Pigeon-12, which is already a substantial challenge for any modern 
solver to date, although the constraint matrix is only $1333 \times 552$
and can be reduced to 738 rows and 288 columns in the presolve.
Modern integer programming solvers fail to solve instances 
larger than this, even considering all the additional constraints described above.

\paragraph{Discretisations}
\label{sec:formulation}

Discretised relaxations proved to be very strong 
in scheduling problems corresponding to one-dimensional packing (\cite{Sousa1992,MR1308272,MR2324961}) and can be shown to be 
asymptotically optimal for various geometric problems both in two dimensions (\cite{MR623065}) and in higher (\cite{MR777998}) dimensions.
\cite{beasley1985exact} has 
extended the formulation to 2D cutting applications,
in the process of deriving a non-linear formulation, for which he proposed solvers.
\cite{Marecek2011} have extended the formulation to the 3D problem of packing boxes
into a larger box, with a considerable amount of work being done independently and subsequently \cite{Junqueira2011,deQueiroz2012,Junqueira2013}.
In this formulation, the small boxes are partitioned into types $t \in \{1,2,\ldots,n\}$, where
boxes of one type share the same triple of dimensions. 
$A_t$ is the number of boxes of type $t$ available.
The large box is discretised into units of space, possibly non-uniformly,
with the indices $(x, y, z) \in D \subset \R^3$ of used to index 
the space-indexed binary variable:

\begin{align}
\label{eq:muvars}
\mu_{x,y,z}^{t}  = \begin{cases}
   \; 1 & \text{if a box of type } t \text{ is placed such that its lower-back-left} \\
        & \text{ vertex is at coordinates } x, y, z \\
   \; 0 & \text{otherwise}  
 \end{cases}
\end{align}

Without allowing for rotations, the formulation reads:

\begin{align}
\max & \sum_{x,y,z,t}{\mu_{xyz}^t w_t} \label{objective-discretised} \\
\text{s.t. } \sum_t{\mu_{xyz}^t} & \leq 1 \quad \forall x,y,z \label{muconstraint1} \\
\mu_{xyz}^t & = 0 \quad \forall x,y,z,t \text{ where } x+L_{1t} > D_X \text{ or } \notag \\
            & \quad \quad \text{ or } y+L_{2t} > D_Y \text{ or } z+L_{3t} > D_Z \label{containedconstraint}\\
\sum_{x,y,z}{\mu_{xyz}^t} & \leq A_t \quad \forall t \label{muconstraint2} \\
\sum_{x',y',z',t' \in f(D, n, L, x, y, z, t)} \mu_{x'y'z'}^{t'} &\leq 1 \quad \forall x,y,z,t \label{no-overlap}\\
\mu_{xyz}^t & \in \{ 0, 1\} \quad \forall x,y,z,t
\end{align}
where one may use Algorithm~\ref{algo:genNonOverlap} or similar\footnote{
For one particularly efficient version of Algorithm~\ref{algo:genNonOverlap}, 
see class {\tt DiscretisedModelSolver} and its method {\tt addPositioningConstraints} in file 
{\tt DiscretisedModelSolver.cpp}
available at \url{http://discretisation.sourceforge.net/spaceindexed.zip} (September 30th, 2014).
This code is distributed under GNU Lesser General Public License.
}
 to generate the index set $f$ in Constraint~\ref{no-overlap}.


The constraints are very natural: No region in space may be occupied by more than one box type (\ref{muconstraint1}),
boxes must be fully contained within the container (\ref{containedconstraint}),
there may not be more than $A_t$ boxes of type $t$ (\ref{muconstraint2}),
and boxes cannot overlap (\ref{no-overlap}). 
There is one non-overlapping constraint (\ref{no-overlap}) for each discretised unit of space and type of box.

\begin{algorithm}[tb]
\caption{
$f(D, n, L, x, y, z, t)$
}
\label{algo:genNonOverlap}
\begin{algorithmic}[1]
  \STATE \textbf{Input:} Discretisation as indices $D \subset \R^3$ of $\mu$ variables \eqref{eq:muvars}, number of boxes $n$, dimensions $L_{\alpha i} \in \R \quad \forall \alpha = x, y, z, i=1, \ldots, n$, 
indices $(x, y, z) \in D$, $t \in \{ 1, 2, \ldots, n\} $ of one scalar within $\mu$ \eqref{eq:muvars} 
  \STATE \textbf{Output:} Set of indices of $\mu$ variables \eqref{eq:muvars} to include in a set-packing inequality \eqref{no-overlap} \rule{0pt}{2.5ex}
  \vspace{1ex}
   \STATE Set $S \leftarrow \{ (x, y, z, t) \} $
   \FOR{ each other unit $(x', y', z') \in D, (x, y, z) \not = (x', y', z')$ } 
   \FOR{ each box type $t' \in \{1, 2, \ldots, n \} $ }  
    \IF{ box of type $t$ at $(x, y, z)$ overlaps box of type $t'$ at $(x', y', z')$, i.e., 
		$(x \le x' + L_{x't'} \le x + L_{xt}) \land (y \le y' + L_{y't'} \le y + L_{yt})$ \\
     \quad \quad $\land (z \le z' + L_{z't'} \le z + L_{zt}) $
} 
        \STATE $S \leftarrow S \cup \{ (x', y', z', t') \}$
    \ENDIF
   \ENDFOR
   \ENDFOR
   \RETURN S
\end{algorithmic}
\end{algorithm}   

In order to support rotations, new box types need to be generated for each allowed rotation and 
linked via set packing constraints, which are similar to Constraint~\ref{muconstraint2}.
In order to extend the formulation to the van loading problem, 
it suffices to add the payload capacity constraint $\sum_{x,y,z,t}{\mu_{xyz}^t m_t} \leq p$.

\section{Finding the Precedence-Constrained Component}

In the rest of the chapter, our goal is to extract the precedence-constrained component
in Chen/Padberg/Fasano formulation from a general integer linear program,
and reformulate it into the discretised formulation.
First, let us state the problem of extracting the Chen/Pad\-berg/Fa\-sa\-no relaxation formally:

\begin{problem}
\label{prob:prec-comp}
{\sc Precedence-Constrained Component/Row-Block Extraction}:
Given positive integers $d, m_1, m_2$, an $m_1 \times d$ integer matrix $A_1$, an $m_2 \times d$ integer matrix $A_2$,  
      an $m_1$-vector $b_1$ of integers, and an $m_2$-vector $b_2$ of integers,  
      corresponding to a mixed integer linear program with constraints $A_1 y = b_1$, $A_2 y \le b_2$, 
find the largest integer $n$ (``the maximum number of boxes''), such that: 
\begin{itemize}  
\item there exists a $4n \times 9n$ submatrix/row-block $E$ of $A_1$, 
               corresponding only to binary variables, which we denote $\delta$,
               with zero coefficients elsewhere in the rows  
\item there exists a $9n(n-1)/2 + 6n + 1 \times 3n(n-1)/2 + 12n$ submatrix/row-block $F$ of $A_2$,
               corresponding to $9n$ binary variables $\delta$ as before,
                                $3n(n-1)/2$ binary variables denoted $\lambda$, and 
                                $3n$ continuous variables denoted $x$, 
               with zero coefficients elsewhere in the rows
\item $E$ contains $n$ rows with exactly 2 non-zero coefficients $\pm 1$, corresponding to (\ref{CP:eq6nzs:1}), and 0 in the right-hand side $b_1$
\item $E$ contains $3n$ rows with exactly 4 non-zero coefficients $\pm 1$, corresponding to (\ref{CP:eq6nzs:2}), and 0 in the right-hand side $b_1$
\item $F$ contains $6n$ rows with exactly 5 non-zero coefficients, some not necessarily $\pm 1$, corresponding to (\ref{CP:5nzs:1}--\ref{CP:5nzs:2}), and 0 in the right-hand side $b_2$
\item $F$ contains $6n(n-1)$ rows with exactly 9 non-zero coefficients, some not necessarily $\pm 1$, corresponding to (\ref{CP:12nzs:1}--\ref{CP:12nzs:2}), and 0 in the right-hand side $b_2$
\item $F$ contains $2n(n-1)$ rows with exactly 9 non-zero coefficients $\pm 1$, corresponding to (\ref{CP:9nzs:1}--\ref{CP:9nzs:2}), and 0 in the right-hand side $b_2$ 
\item $F$ contains $n(n-1)$ rows with exactly 12 non-zero coefficients $\pm 1$, corresponding to (\ref{CP:12nzs}), and 1 in the right-hand side $b_2$
\item $F$ contains $1$ rows with exactly $3n$ non-zero coefficients, not necessarily $\pm 1$, corresponding to (\ref{CP:3nnzs}), and a positive number in the right-hand side $b_2$.
\end{itemize}  
\end{problem}

Notice that by maximising the number of rows involved, we also maximise the number
of boxes, as number $r = 9n(n-1)/2 + 6n + 1$ of rows is determined by number $n$ of boxes.
The following can be seen easily:


\begin{theorem*}
{\sc Precedence-Constrained Row-Block Extraction} is in $\mathcal{P}$.
\end{theorem*}

\begin{proofsketch} 
A polynomial-time algorithm for extracting the pre\-ce\-dence-con\-strain\-ed block can 
clearly rely on there being a polynomial number of blocks of the required size.
See Algorithm~\ref{algo:PrecedenceConstrainedBlock}.
\end{proofsketch}

Algorithm~\ref{algo:PrecedenceConstrainedBlock} displays a very general algorithm
schema for {\sc Precedence-Constrained Row-Block Extraction}. Notably, the test of Line~\ref{pcb:F-test}
requires elaboration.
First, one needs to partition the block into the 5 families of rows (\ref{CP:5nzs:1}--\ref{CP:3nnzs}). 
Some rows (\ref{CP:9nzs:1}--\ref{CP:9nzs:2}, \ref{CP:5nzs:1}--\ref{CP:5nzs:2}, \ref{CP:3nnzs}) are clearly 
determined by the numbers of non-zeros (9, 5, and 3). 
One can distinguish between others (\ref{CP:12nzs:1}--\ref{CP:12nzs:2} and \ref{CP:12nzs}) by their 
right-hand sides.
The test as to whether the rows represent the constraints (\ref{CP:5nzs:1}--\ref{CP:3nnzs}) is based on
identifying the variables.
Continuous variables $x$ are, however, identified easily and binary variables $\delta$ 
are determined in Line~\ref{pcb:E-test}. What remains are variables $\lambda$.

\begin{algorithm}[tb]
\caption{
{\tt PrecedenceConstrainedBlock}($A_1,b_1,A_2,b_2$)
}
\label{algo:PrecedenceConstrainedBlock}
\begin{algorithmic}[1]
  \STATE \textbf{Input:} $A_1 x = 1, A_2 x \le 1$, that is $m_1 \times d$ matrix $A_1$ and $m_2 \times d$ matrix $A_2$, $m_1$-vector $b_1$, $m_2$-vector $b_2$  
  \STATE \textbf{Output:} Integer $k$ and blocks $E, F$ of $A_1, A_2$ \rule{0pt}{2.5ex}
  \vspace{1ex}
\STATE Set $k_{\mbox{max}}$ to the largest integer $k$ such that there are $k$ subsequent rows in $A_1$ with exactly 6 non-zero elements, all $\pm 1$
\FOR{ integer $k = 4n$ from $ k_{\mbox{max}}$ down to 4 } 
 \FOR{ $4n \times 9n$ block $E$ in $A_1$ such that all rows have 6 non-zeros $\pm 1$ } \label{Step1}
  \IF{ there are no $6n$ other rows $A_1$ in corresponding to \eqref{CP:eq6nzs:2} } \label{Step2}
      \label{pcb:E-test}
      \STATE Continue
  \ENDIF
   \FOR{ $7n + 9n(n-1)/2 \times 3n^2 + 9n$ block $F$ in $A_2$ } \label{Step3}
    \IF{ $F$ cannot be partitioned into (\ref{CP:5nzs:1}--\ref{CP:3nnzs}) }    
      \label{pcb:F-test}
      \STATE Continue
    \ENDIF
    \RETURN $n, E, F$
   \ENDFOR
 \ENDFOR                              
\ENDFOR
\end{algorithmic}
\end{algorithm}   

\begin{figure*}[p]
	\centering
	\caption{The workings of Algorithm 1 illustrated on instance Pigeon-02 in the Chen/Padberg/Fasano formulation,
	as introduced in Figure~\ref{pigeon02}.
  }
  \vskip 5mm
  \begin{tabularx}{\textwidth}{ XXX }
	\includegraphics[scale=1.0,trim={0cm 0cm 0cm 0cm},clip]{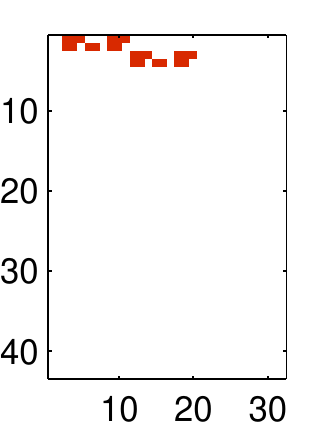} &
	\includegraphics[scale=1.0,trim={0cm 0cm 0cm 0cm},clip]{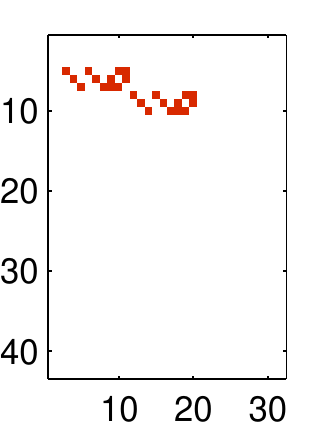} &
	\includegraphics[scale=1.0,trim={0cm 0cm 0cm 0cm},clip]{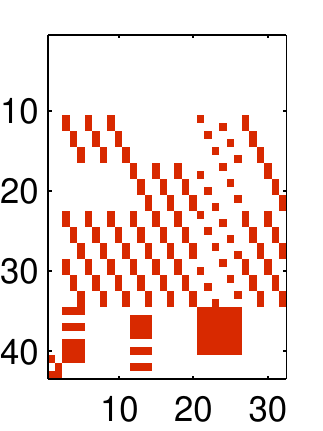}\\
	Line \ref{Step1}: Block $E$ with equalities \eqref{CP:eq6nzs:1} &
  Line \ref{Step2}: The remaining equalities \eqref{CP:eq6nzs:2} &
  Line \ref{Step3}: The remaining inequalities (\ref{CP:5nzs:1}--\ref{CP:12nzs} ) \\[5mm] 
	\includegraphics[scale=1.0,trim={0cm 0cm 0cm 0cm},clip]{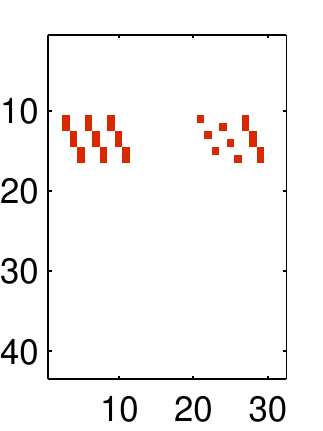} &
	\includegraphics[scale=1.0,trim={0cm 0cm 0cm 0cm},clip]{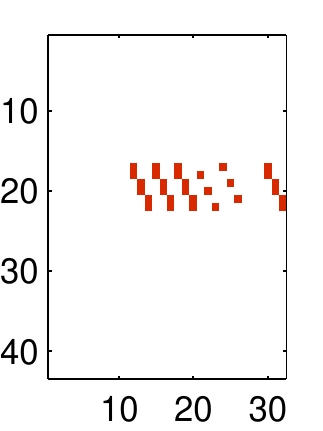} &
	\includegraphics[scale=1.0,trim={0cm 0cm 0cm 0cm},clip]{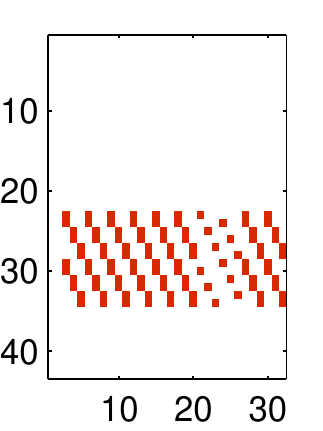}\\
Inequalities \eqref{CP:5nzs:1} &
Inequalities \eqref{CP:5nzs:2} &
Inequalities (\ref{CP:12nzs:1}--\ref{CP:12nzs:2}) \\[5mm]
	\includegraphics[scale=1.0,trim={0cm 0cm 0cm 0cm},clip]{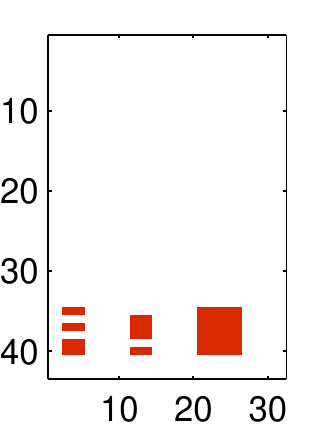} &
	\includegraphics[scale=1.0,trim={0cm 0cm 0cm 0cm},clip]{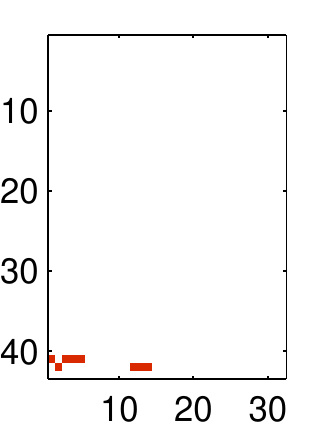} &
	\includegraphics[scale=1.0,trim={0cm 0cm 0cm 0cm},clip]{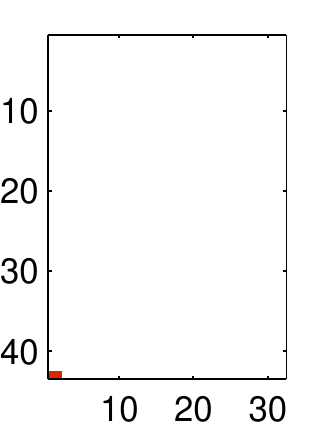}\\
Inequalities (\ref{CP:9nzs:1}--\ref{CP:9nzs:2}) &
Inequalities \eqref{CP:12nzs} &
Inequality \eqref{CP:3nnzs}
	\end{tabularx}
	\label{pigeon02_works}
\end{figure*}

\section{Exploiting the Packing Component}

In order to reduce the number of regions of space, and thus the number of variables 
in the formulation, a sensible space-discretisation method should be employed. 
In many transport applications, for instance, there are only a small number of package types, 
with the ISO 269 standard giving the dimensions of the package. Using the discretisation of one millimeter, one could indeed introduce $162 \times 229 \times h$ variables to represent a package with the C5 base and height $h$ millimeters, but if there are only packages with base-sizes specified by ISO 216 standard and larger than C5 to be packed in the batch, it would make sense to discretise to units of space representing 162 x 229 mm in two dimensions. 
The question is how to derive such a discretisation in a general-purpose system.

The greatest common divisor (GCD) reduction can be applied on a per-axis basis, 
finding the greatest common divisor between the length of the container for an axis 
and all the valid lengths of boxes that can be aligned along that axis and scaling 
by the inverse of the GCD. 
This is trivial to do and is useful when all lengths are multiples of a large number, 
which may be common in certain situations.

In other situations, this may not reduce the number of variables at all, and it
may be worth tackling the optimisation variant of:

\begin{problem}
\label{prob:discretisation}
The {\sc Discretisation Decision}: Given integer $k$, dimensions of a large box (``container'') $x, y, z > 0$ and dimensions of
$n$ small boxes $D \in {\mathbbm{R}}^{n \times 3}$ with associated values $w \in {\mathbbm{R}}^{n}$,
and specification $r = \{0, 1\}^{n \times 6}$ of what rotations are allowed,  
decide if there are $k$ points, where the lower-bottom-left vertex of each box can be positioned in any optimum solution of
{\sc Container Loading Problem}. Formally, $k = |I_X| |I_Y| |I_Z|$, where index-sets $I_X, I_Y, I_Z$ refer to the set 
$\{1c, 2c, \ldots, \max\{x, y, z\} c^{-1} \}$, where $c$ is the greatest common divisor across elements in $D$.  
\end{problem}

Considering the following:

\begin{theorem*}[\cite{MR819147}]
The decision whether the optimum of an instance of {\sc Knapsack} is unique is $\Delta_2^p$-Complete,
where $\Delta_2^p$ is the class of problems that can be solved in polynomial time using oracles from $\mathcal{NP}$.
\end{theorem*} 

\begin{theorem*}
{\sc Discretisation Decision} is $\Delta_2^p$-Hard.
\end{theorem*}
 
\begin{proofsketch}
One could check for the uniqueness of the optimum of an instance of {\sc Knapsack} using any algorithm for {\sc Discretisation Decision} (Problem~\ref{prob:discretisation}). 
\end{proofsketch}
 
We can, however, use non-trivial non-linear space-discretisation heuristics. 
Early examples include \cite{herz1972recursive}.
We use the same values as before on a per-axis basis, i.e. the lengths 
of any box sides that can be aligned along the axis. 
We then use dynamic programming to generate all valid locations for a box to be placed.
See Algorithm~\ref{algo:DiscretisationDP}.  
For example, given an axis of length 10 and box lengths of 3, 4 and 6, we can place 
boxes at positions 0, 3, 4, 6, and 7. 8, 9 and 10 are also possible, 
but no length is small enough to still lie within the container if placed at these points. 
This has reduced the number of regions along that axis from 10 to 5. 
An improvement on this scale may not be particularly common in practice, though the approach can help where the GCD is 1.
It is obvious that this approach can be no worse than the GCD method at discretising the container. 
This also adds some implicit symmetry breaking into the model. 
Notice that the algorithm runs in time polynomial in the size of the output it produces,
  but this is may be exponential in the size of the input and considerably more than
  the size of the best possible output.

\begin{algorithm}[tb]
\caption{{\tt DiscretisationDP}($M, n, D, r$)}
\label{algo:DiscretisationDP}
\begin{algorithmic}[1]
  \STATE \textbf{Input:} Dimensions $M \in {\mathbbm{R}}^3$ of the container, dimensions of $n$ small boxes $D \in {\mathbbm{R}}^{n \times 3}$,
                         and allowed rotations $r = \{0, 1\}^{n \times 6}$    
  \STATE \textbf{Output:} Integer $k$ and $k$ possible positions \rule{0pt}{2.5ex}
  \vspace{1ex}
\STATE $P = \emptyset$
\FOR{ axis with limit $m \in M$ }
 \STATE $L = \{ d \; | \; d \in D $ may appear along this axis, given allowed rotations $ r \} $  
 \STATE $P = $ closure of $ P \cup \{ p + l \; | \; p \in P, l \in L, p + l \le m \} $,\\
        optionally pruning $p \in P$ that cannot occur due to the dependence of the axes and the fact each box can be packed at most once
\ENDFOR
\RETURN $|P|, P$
\end{algorithmic}
\end{algorithm}   

\section{Computational Experience}
\label{sec:computational}

The formulations were tested on two sets of instances, introduced in \cite{Marecek2011}:
\begin{itemize}
\item 3D Pigeon Hole Problem instances, Pigeon-$n$, where $n + 1$ unit cubes are to be packed
into a container of dimensions $(1 + \epsilon) \times (1 + \epsilon) \times n$. 
\item SA and SAX datasets, which are used to test the dependence of solvers' performance 
on parameters of the instances, notably the number of boxes, heterogeneity of the boxes, and physical dimensions of the container.
There is 1 pseudo-randomly generated instance for every combination of container sizes ranging from 5--100 in steps of 5 
units cubed and the number of boxes to pack ranging from 5--100 in steps of 5. 
The SA datasets are perfectly packable, i.e. are guaranteed to be possible to load the container 
with 100\% utilisation with all boxes packed. 
The SAX are similar but have no such guarantees; the summed volume of the boxes is greater than that of the container. 
\end{itemize}
All of the instances are available at \url{http://discretisation.sf.net}.
Some of these instances have been included in MIPLIB 2010 by \cite{KochEtAl2011} and have been widely utilised in 
benchmarking of integer programming solvers ever since.

\begin{table}[tb]
  \centering
  \caption{The performance of various solvers on 3D Pigeon Hole Problem instances encoded in the Chen/Padberg/Fasano formulation. 
  ``--'' denotes that optimality of the incumbent solution has not been proven within an hour.}
  \vskip 5mm
    \begin{tabular}{rrrr}
    \addlinespace
    \toprule
          &       & Time (s) &  \\
    \midrule
          & Gurobi 4.0 & CPLEX 12.4 & SCIP 2.0.1 + CLP \\
    Pigeon-01 & \multicolumn{1}{c}{$<$ 1} & \multicolumn{1}{c}{$<$ 1} & \multicolumn{1}{c}{$<$ 1} \\
    Pigeon-02 & \multicolumn{1}{c}{$<$ 1} & \multicolumn{1}{c}{$<$ 1} & \multicolumn{1}{c}{$<$ 1} \\
    Pigeon-03 & \multicolumn{1}{c}{$<$ 1} & \multicolumn{1}{c}{$<$ 1} & \multicolumn{1}{c}{$<$ 1} \\
    Pigeon-04 & \multicolumn{1}{c}{$<$ 1} & \multicolumn{1}{c}{$<$ 1} & \multicolumn{1}{c}{$<$ 1} \\
    Pigeon-05 & \multicolumn{1}{c}{$<$ 1} & \multicolumn{1}{c}{$<$ 1} & \multicolumn{1}{c}{3.3} \\
    Pigeon-06 & \multicolumn{1}{c}{$<$ 1} & \multicolumn{1}{c}{$<$ 1} & \multicolumn{1}{c}{37.9} \\
    Pigeon-07 & \multicolumn{1}{c}{1.5} & \multicolumn{1}{c}{$<$ 1} & \multicolumn{1}{c}{779.3} \\
    Pigeon-08 & \multicolumn{1}{c}{7.4} & \multicolumn{1}{c}{$<$ 1} & \multicolumn{1}{c}{--} \\
    Pigeon-09 & \multicolumn{1}{c}{88.6} & \multicolumn{1}{c}{66.4} & \multicolumn{1}{c}{--} \\
    Pigeon-10 & \multicolumn{1}{c}{1381.4} & \multicolumn{1}{c}{686.3} & \multicolumn{1}{c}{--} \\
    Pigeon-11 & \multicolumn{1}{c}{--} & \multicolumn{1}{c}{--} & \multicolumn{1}{c}{--} \\
    Pigeon-12 & \multicolumn{1}{c}{--} & \multicolumn{1}{c}{--} & \multicolumn{1}{c}{--} \\
    \bottomrule
    \end{tabular}%

  \label{tab:pigeon-results}
\end{table}

\begin{table}[tb]
  \centering
  \caption{The performance of Gurobi 4.0 on 3D Pigeon Hole Problem instances encoded in the Chen/Padberg/Fasano
  and the discretised formulations as reported in \cite{Marecek2011}. ``--'' denotes that no integer solution has been found. }
  \vskip 5mm
    \begin{tabular}{lrr}
    \addlinespace
    \toprule
          & \multicolumn{2}{c}{Time (s)} \\
    \midrule
          & Chen/Padberg/Fasano & Discretised \\
    Pigeon-01 & \multicolumn{1}{c}{$<$ 1} & \multicolumn{1}{c}{$<$ 1} \\
    Pigeon-02 & \multicolumn{1}{c}{$<$ 1} & \multicolumn{1}{c}{$<$ 1} \\
    Pigeon-03 & \multicolumn{1}{c}{$<$ 1} & \multicolumn{1}{c}{$<$ 1} \\
    Pigeon-04 & \multicolumn{1}{c}{$<$ 1} & \multicolumn{1}{c}{$<$ 1} \\
    Pigeon-05 & \multicolumn{1}{c}{$<$ 1} & \multicolumn{1}{c}{$<$ 1} \\
    Pigeon-06 & \multicolumn{1}{c}{$<$ 1} & \multicolumn{1}{c}{$<$ 1} \\
    Pigeon-07 & \multicolumn{1}{c}{1.5} & \multicolumn{1}{c}{$<$ 1} \\
    Pigeon-08 & \multicolumn{1}{c}{7.4} & \multicolumn{1}{c}{$<$ 1} \\
    Pigeon-09 & \multicolumn{1}{c}{88.6} & \multicolumn{1}{c}{$<$ 1} \\
    Pigeon-10 & \multicolumn{1}{c}{1381.4} & \multicolumn{1}{c}{$<$ 1} \\
    Pigeon-100 & \multicolumn{1}{c}{--} & \multicolumn{1}{c}{$<$ 1} \\
    Pigeon-1000 & \multicolumn{1}{c}{--} & \multicolumn{1}{c}{1.0} \\
    Pigeon-10000 & \multicolumn{1}{c}{--} & \multicolumn{1}{c}{1.8} \\
    Pigeon-100000 & \multicolumn{1}{c}{--} & \multicolumn{1}{c}{4.2} \\
    Pigeon-1000000 & \multicolumn{1}{c}{--} & \multicolumn{1}{c}{45.1} \\
    Pigeon-10000000 & \multicolumn{1}{c}{--} & \multicolumn{1}{c}{664.0} \\
    Pigeon-100000000 & \multicolumn{1}{c}{--} & \multicolumn{1}{c}{--} \\
    \bottomrule
    \end{tabular}%
  \label{tab:pigeontwo}%
\end{table}%


For the 3D Pigeon Hole Problem, results obtained within one hour using three leading solvers and the Chen/Padberg/Fasano formulation 
without any reformulation are shown in 
Table~\ref{tab:pigeon-results}, while Table~\ref{tab:pigeontwo} compares the results on both formulations. 
These tests and further tests reported below were performed on a 64-bit computer running Linux, 
 which was equipped with 2 quad-core processors (Intel Xeon E5472) and 16~GB memory. 
The solvers tested were IBM ILOG CPLEX 12.4, Gurobi Solver 4.0, and SCIP 2.0.1 of \cite{MR2520442} with CLP
as the linear programming solver. 
Pigeon-02 is easy to solve for any modern solver. IBM ILOG CPLEX 12.4 eliminates 17 rows and 26 columns in presolve
and performs a number of further changes. 
The reduced instance has 23 rows and 22 columns and the reported run-time is 0.00 seconds.
Pigeon-10 is the largest instance reliably solvable within an hour, but that should not be surprising, 
considering it has 525 rows, 220 columns, and 3600 nonzeros in the constraint matrix after presolve of CPLEX 12.4.
None of the solvers managed to prove optimality of the incumbent solution for Pigeon-12 within an hour using the Chen/Padberg/Fasano formulation, 
although the instance of linear programming had only 627 rows, 253 columns, and 4268 non-zeros after pre-solve. 
As of September 2014, instances up to pigeon-13 using the Chen/Padberg/Fasano,
have been solved in the process of testing integer programming solvers without the automatic reformulation,
albeit at a great expense of computing time.
In contrast, the reformulation and discretisation makes it possible to solve Pigeon-10000000 within an hour, 
where the instance of linear programming had 10000002 rows, 10000000 columns, and 30000000 non-zeros.
The time for the extraction and reformulation of the instance was under 1 second across of the instances.

\begin{figure*}[p]
	\centering
	\caption{The best solutions obtained within an hour per solver per instance from the SA dataset 
for varying number of boxes (vertical axis) and the length of the side of the container (horizontal axis). 
  The colours highlight the volume utilisation in percent.}
  \vskip -5mm
	\includegraphics[scale=0.7,trim={2cm 5cm 1.5cm 9cm},clip]{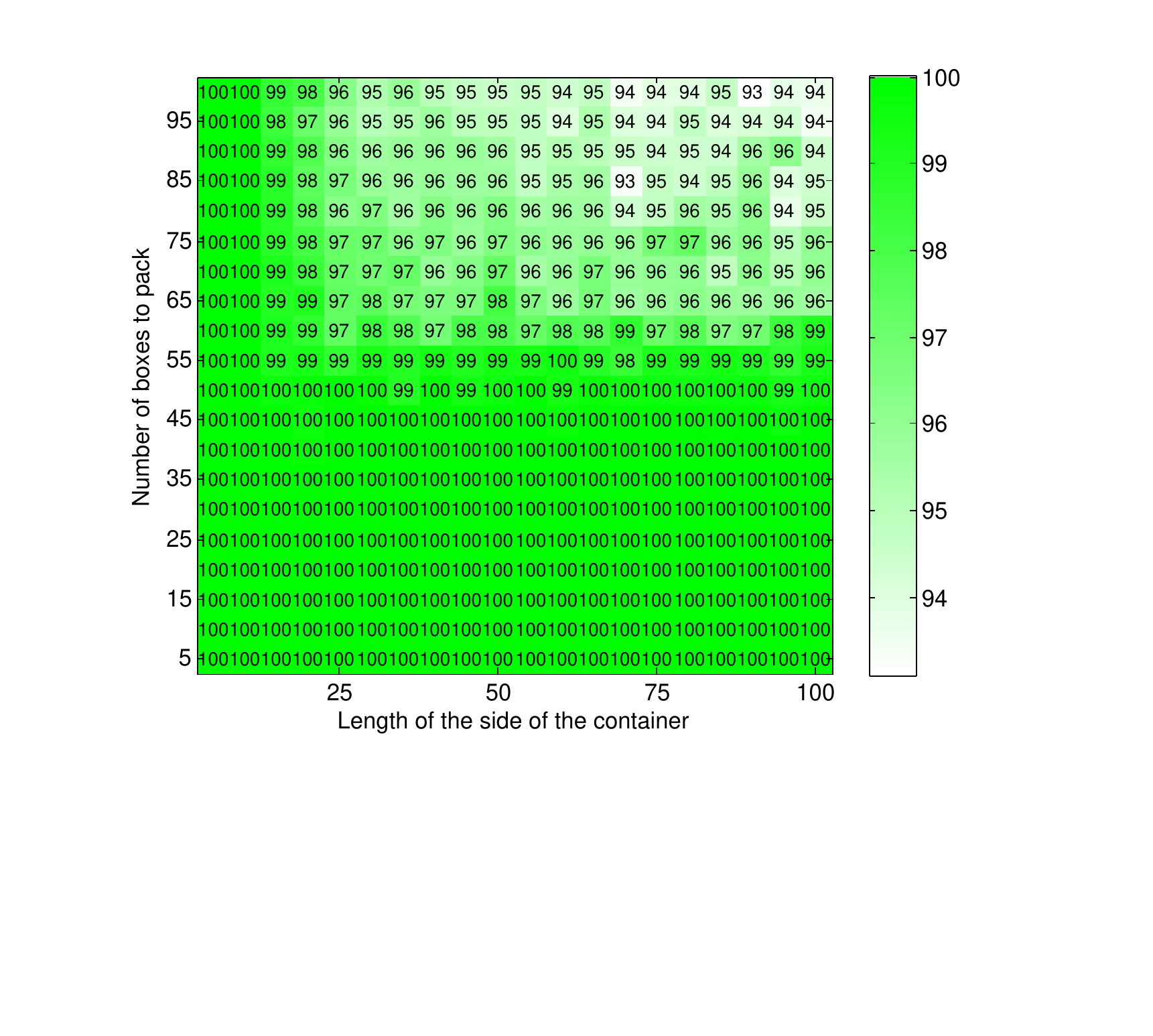}
	\label{contourplot}
  \vskip -30mm
	\caption{As above for the SAX dataset.
 The colours highlight the quality in terms of $100 (1 - s/b)$ for solution with value $s$ and upper bound $b$ after one hour.
  }
  \vskip -5mm
	\includegraphics[scale=0.7,trim={2cm 5cm 1.5cm 9cm},clip]{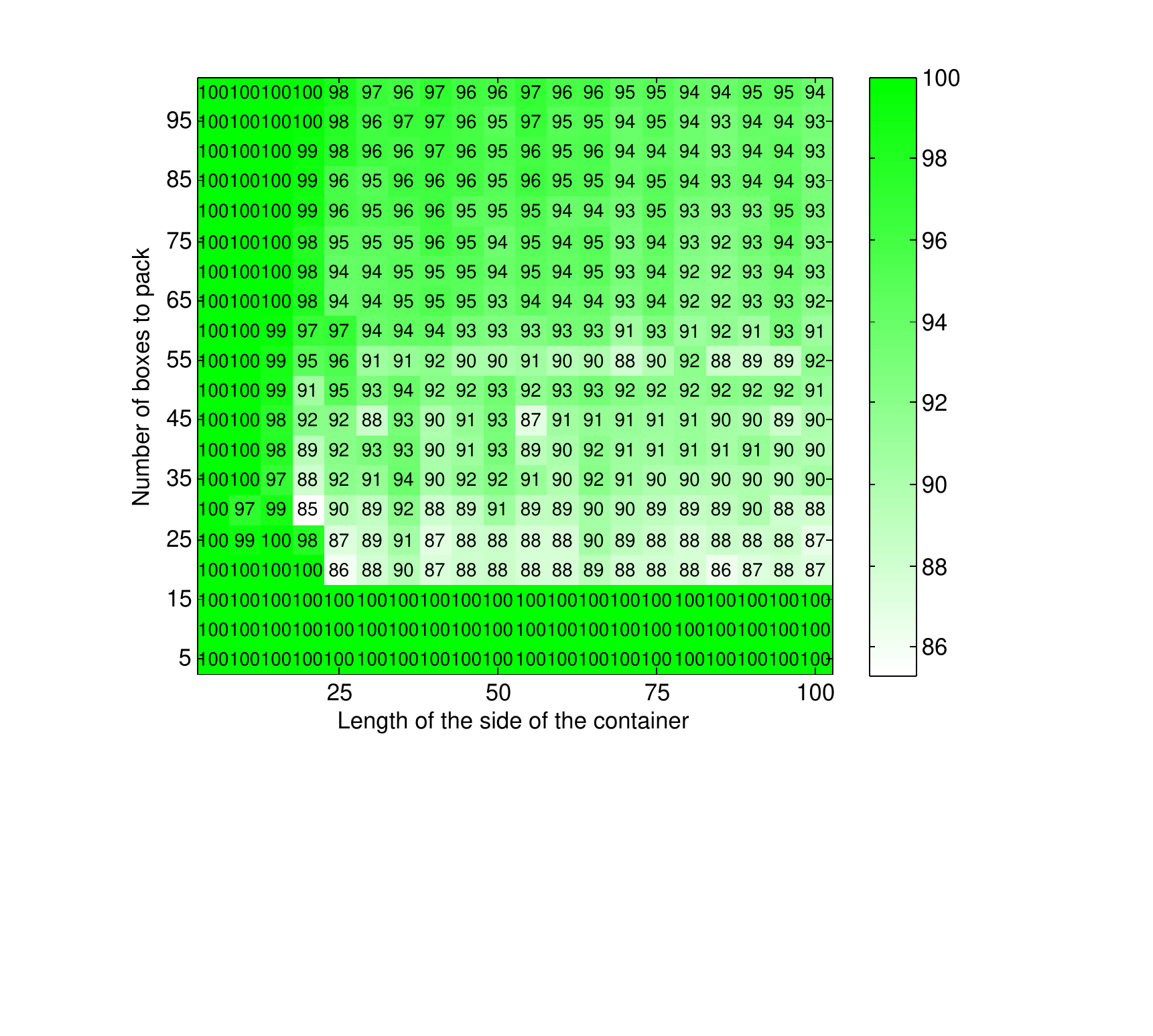}
	\label{contourplotX}
\end{figure*}

For the SA and SAX datasets, Figures \ref{contourplot} and \ref{contourplotX} summarise solutions obtained within an hour using either the Chen/Padberg/Fasano 
formulation or the reformulation, whichever was faster.
This shows that although it may be possible to solve certain instances with 10000000 boxes within an hour to optimality,
real-life instances with hundres of boxes may still be challenging, even considering the reformulation.
Nevertheless, as has been pointed out by \cite{Marecek2011}, the space-indexed relaxation provides a particularly strong upper bound.
The mean integrality gap, or the ratio of the difference between root linear programming relaxation value and optimum to optimum
has been 10.49 \% and 0.37 \% for the Chen/Padberg/Fasano and the space-indexed formulation, respectively, on the SA and SAX instances solved to optimality 
within the time limit of one hour. 

\section{Conclusions}

Overall, the discretisation provides a particularly strong relaxation,
and is easy to reformulate to, programmatically, provided the constraints are a contiguous block of rows, 
as they are, whenever the instances are produced by algebraic modelling languages. 
Clearly, it would be good to develop tests whether the reformulation is worthwhile, as the 
discretised relaxation may become prohibitively large and dense 
for instances with many distinct box-types, and box-types or containers of large sizes in terms of the 
units of discretization.
Alternatively, one may consider multi-level discretisations.
Plausibly, one may also apply similar structure-exploiting approaches to ``components'' other than packing.
\cite{MarecekDiss} studied the graph colouring component, for instance.
This may be open up new areas for research in computational integer programming.


\paragraph*{Acknowledgements}
{\footnotesize
The code for computational testing of the performance of solvers packing problems has code developed by \cite{AllenDiss} for \cite{Marecek2011}
and can be downloaded at \url{http://discretisation.sf.net} (September 30th, 2014).
This material is loosely based upon otherwise unpublished Chapter 8 of the dissertation of \cite{MarecekDiss},
but has benefited greatly from the comments of two anonymous referees,
which have helped the author improve both the contents and the presentation.
The views expressed in this chapter are personal views of the author 
and should not be construed as suggestions as to the product road map of 
IBM products.
}

\bibliographystyle{plainnat}
\bibliography{dr}

\end{document}